\documentclass{kybernetika} 

\usepackage{amsmath,amssymb,mathtools,amsfonts,hyperref,epstopdf} 

\renewcommand\mx{\footnotesize\normalfont}
\renewcommand\mxs{\footnotesize\normalfont\slshape}
\renewcommand\mi{\small\normalfont\itshape}

\usepackage[T1]{fontenc}
\usepackage{lmodern}
\usepackage{microtype}
\usepackage[utf8]{inputenc}

\newtheorem{lemma}{Lemma}[section]
\newtheorem{prop}[lemma]{Proposition}
\newtheorem{coro}[lemma]{Corollary}
\newtheorem{thm}[lemma]{Theorem}

\newtheorem{proposition}[lemma]{Proposition}

\newtheorem{definition}[lemma]{Definition}
\newtheorem{remark}[lemma]{Remark}
\newtheorem{example}[lemma]{Example}

\newcommand\N{\mathbb N}
\newcommand\R{\mathbb R}
\newcommand\Z{\mathbb Z}
\newcommand\A{\mathcal A}
\newcommand\B{\mathcal B}
\newcommand\precalt{\mathrel{\prec_{\mathrm{alt}}}}
\newcommand\succalt{\mathrel{\succ_{\mathrm{alt}}}}
\newcommand\preceqalt{\mathrel{\preceq_{\mathrm{alt}}}}
\newcommand\preclex{\mathrel{\prec_{\mathrm{lex}}}}

\newcommand\IS{{\mathit{IS}}}

\newcommand\Dis{\mathord{\mathcal{D}\mkern-1.3mu\mathit{is}}}

\newcommand\mymax[1]{\mathord{{\max}#1}}
\newcommand\mymin[1]{\mathord{{\min}#1}}
\newcommand\maxAG{\mymax{\A_G}}
\newcommand\minAG{\mymin{\A_G}}
\newcommand\minAL{\mymin{\A_L}}
\newcommand\maxAtwo{\mymax{\A^2}}
\newcommand\minAtwo{\mymin{\A^2}}
\newcommand\maxB{\mymax{\B}}

\DeclarePairedDelimiter{\abs}{\lvert}{\rvert}

\begin{document}

\pagestyle{myheadings}

\title{Greedy and lazy representations\newline in negative base systems}

\author{Tom\'a\v s Hejda, Zuzana Mas\'akov\'a and Edita Pelantov\'a}

\markboth{T. Hejda, Z. Mas\'akov\'a and E. Pelantov\'a}{Greedy and lazy representations in negative base systems}

\contact{Tom\'a\v s}{Hejda}{Doppler Institute for Mathematical Physics and Applied Mathematics
 and Department of Mathematics, FNSPE, Czech Technical University in Prague,
 Trojanova 13, 120\,00 Praha 2. Czech Republic.}{tohecz@gmail.com}
\contact{Zuzana}{Mas\'akov\'a}{Doppler Institute for Mathematical Physics and Applied Mathematics
 and Department of Mathematics, FNSPE, Czech Technical University in Prague,
 Trojanova 13, 120\,00 Praha 2. Czech Republic.}{zuzana.masakova@fjfi.cvut.cz}
\contact{Edita}{Pelantov\'a}{Doppler Institute for Mathematical Physics and Applied Mathematics
 and Department of Mathematics, FNSPE, Czech Technical University in Prague,
 Trojanova 13, 120\,00 Praha 2. Czech Republic.}{edita.pelantova@fjfi.cvut.cz}

\maketitle

\begin{abstract}
We consider positional numeration systems with negative real base $-\beta$,
 where $\beta>1$, and study the extremal representations in these systems,
 called here the greedy and lazy representations.
We give algorithms for determination of minimal and maximal
 $(-\beta)$-representation with respect to the alternate order.
We also show that both extremal representations can be obtained as representations
 in the positive base $\beta^2$ with a non-integer alphabet.
This enables us to characterize digit sequences admissible as greedy and lazy
 $(-\beta)$-representation.
Such a characterization allows us to study the set of uniquely representable numbers.
In the case that $\beta$ is the golden ratio and the Tribonacci constant,
 we give the characterization of digit sequences admissible as greedy and lazy
 $(-\beta)$-representation using a set of forbidden strings.
\end{abstract}

\keywords{numeration systems, lazy representation,
greedy representation, negative base, unique representation}

\classification{11A63, 37B10}

\section{Introduction}

A positional number system is given by a real base $\alpha$ with $\abs{\alpha}>1$ and
 by a finite set of digits $\A\subset \R$, usually called alphabet.
If $x\in \R$ can be expressed in the form $\sum_{i\leqslant k} x_i\alpha^i$
 with coefficients $x_i \in \A$, we say that $x$ has an $\alpha$-representation in $\A$.
The most important assumption on the choice of the alphabet $\A$ is that any
 positive real number $x$ has at least one $\alpha$-representation.
This assumption implies that the cardinality of the alphabet
$\#\A$ is at least $\abs{\alpha}$.

It is well known that if the base $\alpha$ is a positive integer and the alphabet is
 $\A= \{0,1,2, \dots, \alpha-1\}$ then any positive real $x$ has an $\alpha$-representation
 and almost all positive reals (up to a countable number of exceptions) have unique representation.
For example, in the decimal numeration system
\[
    \tfrac12 = 0.5000000\ldots = 0.4999999\ldots
,\qquad\text{whereas}\qquad
    \tfrac13 = 0.333333\ldots
.\]
If the base $\alpha$ is not an integer and the alphabet $\A$ is rich enough to represent
 all positive reals, then almost all $x\geqslant0$ have infinitely many representations and
 one can choose among them ``the nicest'' one from some point of view.%
\footnote{Note that although most people prefer writing $\tfrac12 = 0.5$, the shopkeepers consider
 the representation $0.4999\ldots$ nicer than $0.5000\ldots$ for $x=\frac12$.}

The most studied numeration system with a positive base $\alpha \notin \Z$ uses
 the alphabet $\A=\bigl\{0,1,2,\dots, \lfloor\alpha \rfloor\bigr\}$.
The set of all $\alpha$-representations of $x$ understood as infinite words in the alphabet
 $\A$ is lexicographically ordered.
The lexicographically greatest $\alpha$-representation is considered to be the ``nicest''.
Since this $\alpha$-representation is computed by the so called greedy algorithm,
 it is referred to as the greedy representation.
On the other hand, the lexicographically smallest $\beta$-representation is called lazy.
One might also be interested in other representations in positive base systems, for an extensive
 overview we refer to~\cite{KaSt}.

The study of greedy representations for non-integer bases $\alpha>1$ was initiated by R\'enyi
 in~1957~\cite{Re}.
An interest in lazy representations for bases $\alpha \in (1,2)$ with alphabet $\{0,1\}$ was
 started in 1990 by works of Erd\H{o}s, Jo\'o and Komornik~\cite{ErJoKo}.
A systematic study of lazy representations for all bases $\alpha>1$ can be found in the work
 of Dajani and Kraaikamp \cite{DaKr} from 2002.

Recently, Ito and Sadahiro~\cite{ItSa} introduced a numeration system with a negative base
 $\alpha=-\beta<-1$ and with the alphabet $\A=\bigl\{0,1,2,\dots, \lfloor\beta \rfloor\bigr\}$.
They gave an algorithm for computing a $(-\beta)$-representation
 of $x\in\bigl[\frac{-\beta}{1+\beta}, \frac{1}{1+\beta}\bigr)$ and showed that the natural
 order on $\R$ corresponds to the alternate order on such $(-\beta)$-representations.
Using a negative base, we can represent positive and negative numbers without an additional
 bit for the signum $\pm$.
A family of transformations producing negative base representations of numbers for $1<\beta<2$
 is studied in~\cite{DaKa}.
Among others, it is shown that although none of them gives the maximal representation
 in the alternate order, it is produced by a random algorithm, see Theorem 4.2.\@ in~\cite{DaKa}.
In its proof, one can find out that the greedy representation is obtained by periodic application
 of two transformations.

In this article, we focus on negative bases $-\beta$, $\beta>1$ in general, and deduce analogous
 result without introducing random $(-\beta)$-expansions.
Our main result states that both extremal representations can be obtained using the positive base
 $\beta^2$ and a non-integer alphabet $\B$ by application of a transformation of the form
 $T(x)=\beta^2x-D(x)$, where $D(x)\in\B$
 (Theorems~\ref{t:greedybetanadruhou} and~\ref{t:lazybetanadruhou}).
Note that representations using a non-integer alphabet were considered by Pedicini~\cite{Pedi}.
This enables us to exploit results of~\cite{KaSt} for giving necessary and sufficient conditions
 for identifying sequences admissible as greedy and lazy $(-\beta)$-expansions
 (Theorem~\ref{t:admisgreedy}).
For $\beta=\phi$, the golden ratio, and $\beta=\mu$, the Tribonacci constant, we describe such
 sequences in terms of forbidden strings.
Finally, we illustrate how the characterization of admissible greedy and lazy representations
 can be applied for the study of  uniquely representable numbers (Section~\ref{sec:unique}).

\section{How to obtain \texorpdfstring{$\alpha$}{alpha}-representations of real numbers}

In this chapter we recall a method for finding an $\alpha$-representation of a given number
 with general real base $\alpha$, $\abs{\alpha}>1$.
It is clear that if we are able to find a representation for all numbers $x$ from some bounded
 interval $J \subset \R$ containing $0$, then we are also able to find
 an $\alpha$-representation for any $x \in \bigcup_{k\in\Z} \alpha^k J$,
 i.\,e., for any real number (if $0$ is an interior point of $J$ or the base $\alpha$ is negative)
 or for all positive reals or all negative reals (if $0$ is a boundary point of $J$ and
 the base $\alpha$ is positive).
Our definition below is a restriction of the very general numeration system considered by
 Thurston~\cite{Th}.

\medskip

\begin{definition}\label{JDT}
Given a base $\alpha\in\R$, $\abs{\alpha}>1$, a finite
 alphabet $\A\subset\R$ and a bounded interval $J\ni 0$.
Let $D: J \mapsto \A$ be a mapping such that the transformation defined by
 $T(x) = \alpha x - D(x)$ maps $J$ to $J$.
The corresponding $\alpha$-representation of $x$ is a mapping
 $d=d_{\alpha, J, D}: J\mapsto \A^\N$,
\[
    d_{\alpha, J, D}(x) = x_1x_2x_3x_4 \ldots
,\quad\text{where }
    x_k = D\bigl(T^{k-1}(x)\bigr)
.\]
\end{definition}

Let us comment the previous definition.
By the definition of $T$, for any $x \in J$ we have
 $x = \frac{D(x)}{\alpha} + \frac{T(x)}{\alpha}$.
As the value $T(x)$ is required to belong to $J$ as well, we may use this formula recursively
 and obtain the mentioned $\alpha$-representation of $x$ in the alphabet $\A$:
\begin{equation}\label{rozvoj}
    x = \frac{D(x)}{\alpha} + \frac{D\bigl(T(x)\bigr)}{\alpha^2}
    + \frac{D\bigl(T^2(x)\bigr)}{\alpha^3} + \frac{D\bigl(T^3(x)\bigr)}{\alpha^4} + \dotsb
\end{equation}
Clearly, at any step of this recursion we have
\begin{equation}\label{postupne}
    x = \frac{x_1}{\alpha} + \frac{x_2}{\alpha^2} + \dotsb +
    \frac{x_n}{\alpha^n} + \frac{T^n(x)}{\alpha^n}
.\end{equation}

\begin{example}\label{RenyiGreedy}
In \cite{Re}, R\'enyi defined for $\alpha=\beta >1$ the mappings $T:[0,1)\mapsto[0,1)$
 and $D: [0,1) \mapsto \{0,1,\dots, \lceil \beta \rceil -1\}$ by
\[
    D(x) = \lfloor \beta x\rfloor
\quad\text{and}\quad
    T(x) = \beta x - D(x)
.\]
The $\beta$-representation of $x \in [0,1)$ corresponding to this choice of $D$ and $T$ is
 usually called the greedy $\beta$-expansion of $x$.
\end{example}

More general case of number systems with positive base is studied in~\cite{KaSt},
 where the authors admit non-integer digits and, among other, give condition for
 a digit string to be admissible as a number expansion.
We refer to~\cite{KaSt} for an extensive list of useful literature.

\begin{example}\label{ItoSadahiro}
In \cite{ItSa}, Ito and Sadahiro considered negative bases $\alpha=-\beta$, for any $\beta >1$.
On the interval $J = \bigl[\frac{-\beta} {\beta+1},\frac{1}{\beta+1}\bigr)$ the mappings $T$
 and $D$ are defined as follows:
\[
    D(x) = \bigl\lfloor -\beta x + \tfrac{\beta} {\beta+1}\bigr\rfloor
\quad\text{and}\quad
    T(x) =  -\beta x - D(x)
.\]
The corresponding alphabet, i.\,e., the range of $D$,
 is $\A= \bigl\{0,1,\dots, \lfloor\beta \rfloor\bigr\}$.
Let us mention that for $\beta \notin \N$ the Ito--Sadahiro
alphabet and the R\'enyi
 alphabet coincide.
If that $\beta\in\N$, the digit $\lfloor\beta\rfloor$ occurs only in the infinite
 suffix $\bigl(\lfloor\beta\rfloor\bigr)^\omega$ and every number in $J$ can also
 be represented in base $-\beta$ with the smaller alphabet
 $\A\!\setminus\!\bigl\{\lfloor\beta\rfloor\bigr\}\!=\!\{0,1,\dots,\beta\!-\!1\}$.
However, such a representation cannot be obtained simply using the floor or the ceil function.
\end{example}

The R\'enyi $\beta$-expansion and the Ito--Sadahiro
$(-\beta)$-expansion
 are ``or\-der-pre\-serv\-ing'', 
 provided that we choose a suitable order on the set of representations which they produce.

\begin{definition}\label{def:order}
Let $\A \subset \R$ be a finite alphabet ordered by the natural
order ``$<$'' in $\R$. Let $x_1x_2x_3\ldots$ and $y_1y_2y_3\ldots$
be two different strings
 from $\A^\N$ or $\A^n$ for $n\in\N$, $n\geqslant 1$.
Denote $k = \min\{i\mid x_i \neq y_i\}$.
We write
\begin{itemize}
\item
 $x_1x_2x_3 \ldots \preclex y_1y_2y_3 \ldots$ if $x_k < y_k$ and say that
 $x_1x_2x_3 \ldots $  is smaller  than $y_1y_2y_3 \ldots$ in the lexicographical order;
\item
 $x_1x_2x_3 \ldots \precalt y_1y_2y_3 \ldots$ if $(-1)^k x_k < (-1)^k y_k$
 and say that  $x_1x_2x_3 \ldots $  is smaller than $y_1y_2y_3 \ldots$ in the alternate order.
\end{itemize}
\end{definition}

Let us stress that we always compare only strings of the same length, finite or infinite.

\begin{proposition}\label{usporadani}
Let $\alpha$, $\A$, $J$ and $D$ be as in Definition \ref{JDT}. Let
numbers $x,y\in J$ and let $d_{\alpha, J, D}(x) = x_1x_2x_3
\ldots$
 and $d_{\alpha, J, D}(y)=y_1y_2y_3 \ldots$ be their $\alpha$-representations.
\begin{itemize}
\item If $\alpha > 1$ and $D(x)$ is non-decreasing then
 \[
    x< y
 \quad\iff\quad
    x_1x_2x_3 \ldots \preclex y_1y_2y_3\ldots
 .\]
\item If $\alpha <-1 $ and $D(x)$ is non-increasing then
 \[
    x<y
 \quad\iff\quad
    x_1x_2x_3 \ldots \precalt y_1y_2y_3\ldots
 .\]

\end{itemize}
\end{proposition}

\begin{Proof}
We prove the statement for a negative base $\alpha <-1$.
Let us denote $k =  \min\{ i \mid x_i \neq y_i\}$.
As $x_\ell = y_\ell$ for all $\ell = 1,2,\dots, k-1$ we have according to \eqref{postupne} that
\begin{equation}\label{eq:x1}
    x<y \quad\iff\quad\frac{T^{k-1}(x)}{\alpha^{k-1}} < \frac{T^{k-1}(y)}{\alpha^{k-1}}.
\end{equation}
Since $D$ is non-increasing and $x_k=D\bigl(T^{k-1}(x)\bigr)\neq y_k = D\bigl(T^{k-1}(y)\bigr)$,
 we have
\begin{equation}\label{eq:x2}
    x_k>y_k \quad\iff\quad\frac{T^{k-1}(x)}{\alpha^{k-1}} < \frac{T^{k-1}(y)}{\alpha^{k-1}}.
\end{equation}
First we discuss the case when $k$ is even. Since $\alpha < 0$,
 combining~\eqref{eq:x1} and~\eqref{eq:x2}, we obtain
\[
    x<y \quad\iff\quad x_k<y_k \quad\iff\quad x_1x_2x_3 \ldots \precalt y_1y_2y_3\ldots.
\]
If $k$ is odd, then we obtain
\[
    x<y \quad\iff\quad x_k>y_k \quad\iff\quad x_1x_2x_3 \ldots \precalt y_1y_2y_3\ldots.
\]
The proof of the statement for a positive base is analogous.
\end{Proof}

\bigskip

The above proposition allows one to give a condition of admissibility for a digit string to be
 a representation of a number according to the algorithm presented in Definition~\ref{JDT}.
Such an admissibility condition was given in~\cite{Parry} for the R\'enyi number system
 from Example~\ref{RenyiGreedy}, in~\cite{KaSt} for more general positive base systems
 and in~\cite{ItSa} for the number system from Example~\ref{ItoSadahiro}.
Let us stress that the representations and consequently admissibility conditions depend on the
 given transformation $T$ in the scheme $J,D,T$.

From now on, we focus on a different question.
We do not fix any transformation, but among all representations of a number $x$ with a given
 alphabet, we search for an extremal one.

Given a base $\alpha$ and an alphabet $\A \subset \R$ we put
\[
    J_{\alpha, \A} =
    \Bigl\{\begingroup\textstyle\sum\limits_{i=1}^{\infty}\endgroup a_i \alpha^{-i} \,\Big\vert\,a_i \in \A\Bigr\},
\]
the set of numbers representable with negative powers of $\alpha$ and the alphabet $\A$, and for
 any $x \in J_{\alpha,\A}$, we denote the set of its $\alpha$-representations in $\A$ by
\[
    R_{\alpha, \A}(x) = \Bigl\{ x_1x_2x_3 \ldots  \,\Big\vert\,
    x=\begingroup\textstyle \sum\limits_{i=1}^{\infty}\endgroup x_i\alpha^{-i} \text{ and } x_i \in \A\Bigr\}
.\]

Proposition \ref{usporadani} suggests how to choose a suitable ordering on the
 set $R_{\alpha, \A}(x)$.

\begin{definition}
Let $\alpha$ be a real base with $|\alpha|>1$ and let $\A \subset \R$ be a finite alphabet.
\begin{itemize}
\item
 Let $\alpha <-1$ and $x \in J_{\alpha, \A}$.
 Then the maximal  and minimal elements of $R_{\alpha, \A}(x)$ with respect to the alternate
  order are called the greedy and lazy $\alpha$-representations of $x$ in the alphabet
  $\A$, respectively.
\item
 Let $\alpha >1$ and $x \in J_{\alpha, \A}$.
 Then the maximal and minimal elements of $R_{\alpha, \A}(x)$ with respect to the
  lexicographical order are called the greedy and lazy $\alpha$-representations of $x$ in the alphabet
  $\A$, respectively.
\end{itemize}
\end{definition}

In order to justify that the maximal and minimal elements in $R_{\alpha, \A}(x)$ exist, realize that
 $R_{\alpha, \A}(x)$ is a compact subspace of $\A^{\N}$. For, it is a pre-image of a closed set
 $\{x\}\subset\R$ under the map $x_1x_2x_3\ldots \mapsto \frac{x_1}{\alpha^1}+\frac{x_2}{\alpha^2}+
 \frac{x_3}{\alpha^3}+\dotsb$ which is continuous with respect to the Cantor metric on $\A^{\N}$.

\bigskip

For positive bases $\alpha=\beta>1$, the greedy $\beta$-representation in
 the alphabet $\A=\bigl\{0,1,\dots,\lceil\beta\rceil-1\bigr\}$ corresponds to the $\beta$-representation
 produced by the  R\'enyi greedy algorithm described in Example \ref{RenyiGreedy}.
The method for obtaining the lazy representation in the same alphabet was described in \cite{DaKr}.

\begin{figure}
\centering 
    \raisebox{-\height}{\includegraphics[page=2,scale=0.9]{plots.pdf}}%
\hfil
    \raisebox{-\height}{\includegraphics[page=1,scale=0.9]{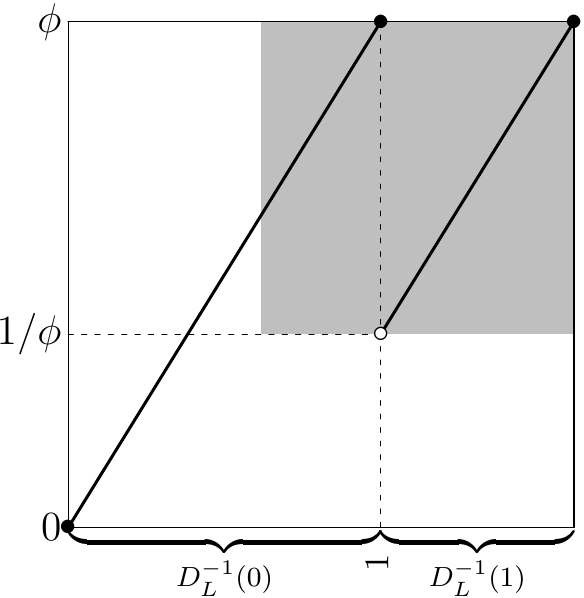}}%
\caption{The greedy and lazy transformations $T_G$ and $T_L$ for $\beta=+\phi$
 with pre-images of $0$ and $1$ under $D_G$ and $D_L$, respectively.}
\label{graf:LGpos}
\end{figure}

Let us demonstrate how to obtain these extremal representations for the case of the golden
 mean $\phi = \frac{1+\sqrt{5}}{2}$.
Since $\lceil\phi\rceil = 2$, we consider the alphabet $\A=\{ 0,1\}$.
As $\phi^2 = \phi + 1$,  we have $\sum_{i=1}^{\infty} \phi^{-i} = \phi$ and so clearly,
\[
    J_{\phi,\A}
    = \biggl\{\sum_{i=1}^{\infty} a_i \phi^{-i} \,\bigg\vert\, a_i \in \{0,1\}\biggr\}
    =[0,\phi]
.\]
If the first digit in a $\phi$-representation of an $x \in J_{\phi, \A}$ is $x_1=0$ then
 necessarily $x \in \frac{1}{\phi} J_{\phi, \A} = [0,1]$.
If the first digit in a $\phi$-representation of $x$ is $x_1=1$ then
 necessarily $x \in \frac{1}{\phi}+ \frac{1}{\phi} J_{\phi, \A}
 = \bigl[\frac{1}{\phi} ,\phi\bigr]$.
It means that
\begin{itemize}
\item
 for $x \in [0, \frac{1}{\phi})$, all representations in $R_{\phi,\A}(x)$ have
 the first digit $x_1=0$;
\item
 for  $x \in ( 1, \phi]$, all representations in $R_{\phi,\A}(x)$ have
 the first digit $x_1=1$;
\item
 for $x \in [ \frac{1}{\phi} ,1]$, both $0$ and $1$ appear as $x_1$.
\end{itemize}
As we now consider the lexicographical order, the greedy representation of
 an  $x \in \bigl[ \frac{1}{\phi} ,1 \bigr]$ has the first digit $x_1=1$ and the lazy
 representation must have $x_1=0$.
Thus the corresponding digit assigning mappings are
\[
    D_G(x) = \begin{cases}
        0& \text{for } x \in \bigl[0,\tfrac{1}{\phi}\bigr)
    ,\\
        1& \text{for } x \in \bigl[\tfrac{1}{\phi}, \phi\bigr]
    ,\end{cases}
\quad\text{and}\quad
    D_L(x) = \begin{cases}
        0& \text{for } x \in [0,1]
    ,\\
        1& \text{for } x \in (1, \phi]
    ,\end{cases}
\]
or equivalently
\[
    D_G(x)=\bigl\lfloor x + \tfrac{1}{\phi^2} \bigr\rfloor
,\qquad
    D_L(x)=\lceil x\rceil -1
\qquad \text{for any }
    x \in (0,\phi)
\]
and
\[
    D_G(0)=D_L(0)=0 ,\qquad D_G(\phi)=D_L(\phi)=1
.\]
The graphs of
the corresponding transformations $T_G$ and $T_L$ are depicted in
Figure~\ref{graf:LGpos}.

\medskip

The transformations $T_G$ and $T_L$ map $[0,\phi] \mapsto [0,\phi]$.
The points $0$ and $\phi$ are fixed points of both transformations.
However, we see that
\begin{itemize}
\item\rule{-2mm}{0mm}
 for any $x \in [0, \phi)$ there exists $n\geqslant 1$ such that
 $T_G^n(x) \in J_G = [0, 1)$ and $T_G(J_G) = J_G$;
\item\rule{-2mm}{0mm}
 for any $x \in (0, \phi]$ there exists $n\geqslant 1$ such that
 $T_L^n(x) \in J_L = (\frac{1}{\phi}, \phi]$ and $T_L(J_L) = J_L$.
\end{itemize}
Moreover the intervals $[0,1)$ and $(\frac{1}{\phi}, \phi]$ are the smallest intervals with
 the above described property.
Such interval is called the attractor of a mapping in \cite{DaKr}.
Let us point out that in both cases there exists an isolated fixed point $\phi=T_G(\phi)$ and
 $0=T_L(0)$ outside the attractors.
In fact, to find the greedy $\beta$-representation, R\'enyi used $J_G=[0,1)$
 as the domain of $D$ and $T$ in the scheme $J,D,T$.
Even when the scheme $J,D,T$ for finding the lazy representation is described,
 the restriction of the mappings $D$ and $T$ to the `attractor' interval of $T$ is preferred.

\section{Extremal representations in negative base systems}

Let us now fix a base $\alpha=-\beta$ for some $\beta>1$, $\beta\notin\N$, and an
 alphabet $\A=\bigl\{0,1,\dots,\lfloor\beta\rfloor\bigr\}$.
(We exclude integer bases because of the phenomena explained in Example~\ref{ItoSadahiro}.)
Using the same arguments as in~\cite{Pedi} it can be shown that the set $I$ of numbers
 representable in this system is an interval, namely
\begin{equation}\label{eq:I}
    I = \Bigl[\frac{-\beta\lfloor\beta\rfloor}{\beta^2-1},
    \frac{\lfloor\beta\rfloor}{\beta^2-1}\Bigr] \eqqcolon [l,r].
\end{equation}

Though the interval $I$ (i.\,e., its boundary points) depends on
the base $\beta$,
 we omit the symbol $\beta$ it in the notation for simplicity.
We denote by $I_a$ the set of numbers which have a $(-\beta)$-representation starting
 with the digit $a\in\A$.
Then
\[
    I_a = \frac{a}{-\beta} + \frac1{-\beta} I=\Bigl[l+\frac{\lfloor\beta\rfloor-a}{\beta},
    r-\frac{a}{\beta}\Bigr]
    = \Bigl[\frac{a}{-\beta}+\frac{r}{-\beta},\frac{a}{-\beta}+\frac{l}{-\beta}\Bigr]
\]
 and $I$ can be written as a (not necessarily disjoint) union of intervals $I=\bigcup_{a\in\A}I_a$.
Obviously, we have $-\beta x - a \in I$  for every $x\in I_a$.

Note that intervals $I_a$ overlap, but only two at a time.
If $\lfloor\beta\rfloor=1$, then it is clear, as they are only two.
If $\lfloor\beta\rfloor\geqslant 2$, then the right endpoint of $I_{a+1}$ is smaller than
 the left endpoint of $I_{a-1}$, namely
\[
    \frac{a+1}{-\beta} + \frac{l}{-\beta}=\frac{a+1}{-\beta} +
    \frac{\lfloor\beta\rfloor}{\beta^2-1} < \frac{a-1}{-\beta} +
    \frac{1}{-\beta}\frac{\lfloor\beta\rfloor}{\beta^2-1}=\frac{a-1}{-\beta}
    + \frac{r}{-\beta}
,\]
which is equivalent to $2\beta-\lfloor\beta\rfloor-2>0$.

Our aim is to provide an algorithm for finding extremal digit strings (with respect to the
 alternate order) representing a given $x$ in the base $(-\beta)$.
Assigning of the first digit of such an extremal representation must prefer --- among
 all possible digits --- the minimal or the maximal one in the alternate order.

We define
\[
    D_m(x)=\begin{cases}
        \lfloor\beta\rfloor & \text{for } x \in I_{\lfloor\beta\rfloor}
    ,\\
        a & \text{for } x \in  I_{a}\setminus I_{a+1}
    ,\ a\in\A,\ a\neq \lfloor\beta\rfloor ,
    \end{cases}\vspace*{-2mm}
\]
and\vspace*{-2mm}
\[
    D_v(x)=\begin{cases}
        0 & \text{for } x \in I_{0}
    ,\\
        a & \text{for } x \in  I_{a}\setminus I_{a-1}
    ,\ a\in\A,\ a\neq 0 ,
    \end{cases}
\]
and the corresponding transformations
\[
    T_m(x) = -\beta x - D_m(x)
    \qquad\text{and}\qquad
    T_v(x) = -\beta x - D_v(x).
\]

It is readily seen that $D_m(x)$ prefers small digits (in the
alternate order) and $D_v(x)$ prefers great digits, i.\,e.,
\begin{equation}
    \label{eq:ven1}
    D_m(x)=a \quad\implies\quad -\beta x-b\notin I\quad\text{for any } b\precalt a,\ b\in\A,
    \vspace*{-2mm}\end{equation}
and\vspace*{-2mm}
\begin{equation}
    \label{eq:ven2}
    D_v(x)=a \quad\implies\quad -\beta x-b\notin I\quad\text{for any } b\succalt a,\ b\in\A.
\end{equation}

Let us stress that in accordance with Definition~\ref{def:order}, alternate order on single digits
 is inverse to the usual one.

\begin{figure}
\centering
    \includegraphics[page=6,scale=1.1]{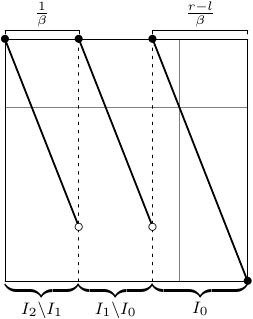}%
\hfil
    \includegraphics[page=7,scale=1.1]{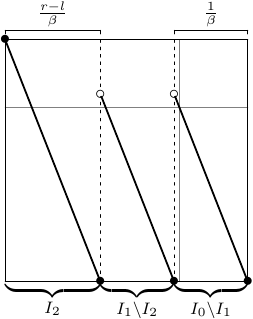}%
\caption{Transformations $T_v$ and $T_m$ which prefer small and large digits in the alternate
 order, respectively (for a base $-\beta\in(-3,-2)$).
The fractions $\frac{1}{\beta}$ and $\frac{r-l}{\beta}$ in the upper part of the figure represent
 two different lengths of intervals indicated in the lower part of the figure.}
\label{graf:malavelka}
\end{figure}

\begin{proposition}\label{prop:algLGobecnebeta}
Let $x \in I$.
\begin{itemize}
\item
 Denote $\varepsilon_0 = x$ and for all $i\geqslant 0$ put
 \[
    z_{2i+1} = D_v(\varepsilon_{2i}),
    \ \varepsilon_{2i+1} = T_v(\varepsilon_{2i})
 \quad\text{and}\quad
    z_{2i+2} = D_m(\varepsilon_{2i+1}),
    \ \varepsilon_{2i+2} = T_m(\varepsilon_{2i+1}).
 \]
 Then $z_1z_2z_3 \ldots$ is the greedy $(-\beta)$-representation of $x$.
\item
 Denote $\eta_0 = x$ and for all $i\geqslant 0$ put
 \[
    y_{2i+1} = D_m(\eta_{2i}),
    \ \eta_{2i+1} = T_m(\eta_{2i})
 \quad\text{and}\quad
    y_{2i+2} = D_v(\eta_{2i+1}),
    \ \eta_{2i+2} = T_v(\eta_{2i+1}).
 \]
 Then $y_1y_2y_3 \ldots$ is the lazy $(-\beta)$-representation of $x$.
\end{itemize}
\end{proposition}

\begin{Proof}
Suppose for contradiction that there exists a representation
$x_1x_2x_3\ldots \in R_{-\beta, \A}(x)$ such that $z_1z_2z_3\ldots
\precalt x_1x_2x_3\ldots$. Let us denote $k=\min\{i \mid z_i \neq
x_i\}$.

If $k=2\ell$ is even, then the inequality $z_1z_2z_3\ldots\precalt
x_1x_2x_3\ldots$
 implies $z_{2\ell} < x_{2\ell}$ and thus $z_{2\ell}=a$ and $x_{2\ell}=b$
 for some $a,b\in{\mathcal A}$, $a<b$.
Since
\begin{multline*}
    \varepsilon_{2\ell-1}
    = (-\beta)^{2\ell-1} \biggl(x - \sum_{i=1}^{2\ell-1} z_i (-\beta)^{-i}\biggr)
    = (-\beta)^{2\ell-1} \biggl(x - \sum_{i=1}^{2\ell-1} x_i (-\beta)^{-i}\biggr)
\\
    = (-\beta)^{2\ell-1} \biggl(x_{2\ell} (-\beta)^{-2\ell}
    + \sum_{i=2\ell+1}^{\infty} x_i (-\beta)^{-i}\biggr)
    = -\frac{1}{\beta} \biggl(b+ \sum_{i=2\ell+1}^{\infty} x_i (-\beta)^{2\ell-i} \biggr)
,\end{multline*} we have $-\beta
\varepsilon_{2\ell-1}\!-\!b\!\in\! I\!=\![l,r]$. On the other
hand, according to \eqref{eq:ven1},
$z_{2\ell}\!=\!D_m(\varepsilon_{2\ell-1})\!=\!a$
 implies that $-\beta \varepsilon_{2\ell-1}- b \notin  I$, which is a contradiction.
The case $k$ odd is analogous.

Very similarly one proves the statement for the lazy representation.
\end{Proof}

\bigskip

Note that for $\beta\in(1,2)$, the description of greedy $(-\beta)$-representations
 of Proposition~\ref{prop:algLGobecnebeta} can be deduced from the proof
 of Theorem~4.2 in~\cite{DaKa}.

\section{Representations in base \texorpdfstring{$\beta^2$}{beta\string^2}
 with non-integer alphabets}

The algorithm for obtaining extremal $(-\beta)$-representations of a number $x$,
 described in Proposition~\ref{prop:algLGobecnebeta}, does not fit in the scheme of
 the Definition~\ref{JDT} for negative base $\alpha=-\beta$.
In particular, there is no transformation $T(x)=-\beta x-D(x)$ which generates for
 every $x$ the greedy (or lazy) $(-\beta)$-representation.
This fact complicates the description of digit strings occurring as greedy or lazy representations.
 Nevertheless, we manage to overcome this handicap.

Defining $T_G\coloneqq T_mT_v$ and $T_L\coloneqq T_vT_m$ we obtain transformations $I\to I$ which
 produce the greedy and lazy $(-\beta)$-representations.
The price to be paid is that the digit assigning functions $D_G$ and $D_L$ are not integer-valued.
For, we have
\begin{gather}\nonumber
    T_G(x) = T_m\bigl(T_v(x)\bigr) = \beta^2 x + \beta D_v(x) - D_m\bigl(-\beta x - D_v(x)\bigr),
\\\shortintertext{i.\,e.,}\label{eq:TG}
    T_G(x) = \beta^2x-D_G(x),
\\\shortintertext{where}\nonumber
    D_G(x) = -\beta D_v(x) + D_m\bigl(-\beta x - D_v(x)\bigr) \in \B,
\end{gather}
 with the alphabet $\B= \{-b\beta+a \mid a,b\in\A\}$, which has $(\#\A)^2$ distinct
 elements, as $\beta\notin\N$.

Let us describe the mappings $T_G$ and $D_G$ (resp. $T_L$ and $D_L$) explicitly.
As $T_G$ is a composition of two piecewise linear functions, $T_G$ is also of this kind.
Its points of discontinuity are described as follows,
\begin{equation}\label{eq:discGpredbezne}
    \Dis_G = \Dis_v \cup T_v^{-1}\bigl(\Dis_m\bigr),
\end{equation}
 where $\Dis_v$ and $\Dis_m$ are points of discontinuity of $T_v$ and $T_m$ respectively.
Clearly, $\#\Dis_v = \#\Dis_m = \#\A- 1$.
We can easily see that
\begin{equation}\label{eq:discm}
    \Dis_m = \Bigl\{\frac{a}{-\beta} + \frac{l}{-\beta}\,\Big\vert\,a\in\A,\ a\neq 0\Bigr\}
\end{equation}
 and
\begin{equation}\label{eq:discv}
    \Dis_v = \Bigl\{\frac{a}{-\beta} + \frac{r}{-\beta}
    \,\Big\vert\,a\in\A,\ a\neq\lfloor\beta\rfloor\Bigr\}.
\end{equation}
As every discontinuity point $\frac{a}{-\beta} + \frac{l}{-\beta}$ of $\Dis_m$ lies in $(r-1,r]$
 and for every letter $a\in\A$ one has $(r-1,r]\subset T_v(I_a)$, the set $T^{-1}_v(\Dis_m)$
 has $\#\A\times (\#\A-1)$ points.
Thus the union in~\eqref{eq:discGpredbezne} has $(\#\A)^2-1$ elements.
The interval $I$ therefore splits into $(\#\A)^2$ subintervals on which the transformation $T_G$
 is continuous.
Explicitly, for every $a\in\A$, $a\neq 0$, we have
\[
    T_v^{-1}\Bigl(\frac{a}{-\beta} + \frac{l}{-\beta}\Bigr) =
    \Bigl\{\frac{b}{-\beta} + \frac{a}{(-\beta)^2}+\frac{l}{(-\beta)^2}
    \,\Big\vert\,b\in\A\Bigr\}.
\]
As $-\beta r=l$, every discontinuity point $\frac{b}{-\beta} + \frac{r}{-\beta}$ of $\Dis_v$
 can be written in the form  $\frac{b}{-\beta} + \frac{l}{(-\beta)^2}$, and thus
\[
    \Dis_G = \Bigl\{\frac{b}{-\beta} + \frac{a}{(-\beta)^2}+\frac{l}{(-\beta)^2}
    \,\Big\vert\,a,b\in\A,\ (b,a)\neq(\lfloor\beta\rfloor,0)\Bigr\}.
\]
In order to obtain an explicit description of the mapping $T_G$, let us define
\begin{equation}\label{eq:gamma}
    \gamma_{ba}\coloneqq l+\frac{\lfloor\beta\rfloor-b}{\beta} = \frac{b}{-\beta} +
    \frac{a}{(-\beta)^2}+\frac{l}{(-\beta)^2},
    \quad\text{for any }a,b\in\A.
\end{equation}
Note that if $(b,a)=(\lfloor\beta\rfloor,0)$, then $\gamma_{ba}=l$.
The set of words of length two over the alphabet $\A$ can be ordered by the alternate order, namely
\begin{multline*}
    \lfloor\beta\rfloor0
    \precalt\lfloor\beta\rfloor1
    \precalt\cdots\precalt
    \lfloor\beta\rfloor\lfloor\beta\rfloor
    \precalt(\lfloor\beta\rfloor-1)0
\precalt\cdots\\\cdots
    \precalt1\lfloor\beta\rfloor
    \precalt00\precalt\cdots
    \precalt0\lfloor\beta\rfloor.
\end{multline*}
Taking into account that $\beta\notin\N$, one can easily see that
\[
    ab\precalt cd \quad\iff\quad \gamma_{ab} = \frac{a}{-\beta}+
    \frac{b+l}{(-\beta)^2} < \frac{c}{-\beta}+ \frac{d+l}{(-\beta)^2}
    = \gamma_{cd}.
\]
In particular, this means that the assignment $ab\leftrightarrow-a\beta+b$,
 for $a,b\in\A$, is a bijection.
On the set $\A^2$ we can define the predecessor and successor functions by
\begin{align*}
    P(ab)&\coloneqq \max\{cd \mid cd\in\A^2,\ cd\precalt ab\}
,\\
    S(ab)&\coloneqq \min\{cd \mid cd\in\A^2,\ cd\succalt ab\}
,\end{align*}
 where maximum and minimum is taken again with respect to the alternate order.
Of course, the predecessor is not defined for the minimal
 element $\minAtwo=\lfloor\beta\rfloor0$ and the successor is not defined for the maximal
 element $\maxAtwo=0\lfloor\beta\rfloor$.

Now we are ready to provide an explicit description of the digit assigning map
 $D_G: I \to \B=\{-b\beta+a\mid a,b\in\A\}$, namely
\begin{equation}\label{eq:DG}
    D_G(x) = \begin{cases}
        -b\beta+a,&\text{if }x\in\bigl[\gamma_{ba},\gamma_{S(ba)}\bigr)
        \text{ for some } ba\neq 0\lfloor\beta\rfloor
    ,\\
        \quad\lfloor\beta\rfloor,&\text{if }x\in\bigl[\gamma_{ba},r\bigr]
        \text{ for } ba= 0\lfloor\beta\rfloor
    .\end{cases}
\end{equation}

Define a morphism $\psi:\B^*\to\A^*$ by
\[
    \psi(-b\beta+a) = ba.
\]
Note that $\psi$ is well defined as $\beta\notin\N$.
The action of the morphism $\psi$ can be naturally extended to infinite words by
 setting $\psi(x_1x_2x_3\ldots) = \psi(x_1)\psi(x_2)\psi(x_3)\ldots$.
All the above considerations can be summarized in the following theorem.

\begin{thm}\label{t:greedybetanadruhou}
Let $\beta>1$, $\beta\notin\N$, $\A=\bigl\{0,1,\dots,\lfloor\beta\rfloor\bigr\}$.
For an $x\in I$, denote by $d_G(x)=d_{\beta^2,I,D_G}(x)$ the $\beta^2$-representation
 of $x$ obtained in accordance with Definition~\ref{JDT} by the
 transformation $T_G:I\to I$, defined by~\eqref{eq:TG}.
Then
\begin{itemize}
\item
 $d_{G}(x)$ is the greedy $\beta^2$-representation of $x$ in the alphabet $\B$;
\item
 $\psi\bigl(d_{G}(x)\bigr)$ is the greedy $(-\beta)$-representation of $x$ in the alphabet $\A$.
\end{itemize}
\end{thm}
\medskip

For the sake of completeness, let us provide an explicit description of the lazy
 transformation $T_L=T_v T_m : I\to I$ and the corresponding digit assigning
 map $D_L : I\to \B$, such that
\begin{equation}\label{eq:TL}
    T_L(x)=\beta^2x-D_L(x),
\end{equation}
 for $x\in I$.
Defining
\[
    \delta_{ba}\coloneqq\frac{b}{-\beta} + \frac{a}{(-\beta)^2}+\frac{r}{(-\beta)^2} ,
    \quad\text{for any }a,b\in\A.
\]
we can write the set of discontinuity points of $T_L$ as
\[
    \Dis_L = \bigl\{\delta_{ba}\,\big\vert\,a,b\in\A,\ ba\neq0\lfloor\beta\rfloor\bigr\}.
\]
Note that for $ba=0\lfloor\beta\rfloor$ one has $\delta_{ba}=r$.
Again,
\[
    ab\precalt cd \quad\iff\quad \delta_{ab} < \delta_{cd}.
\]
We now have the digit assigning function $D_L: I \to \B=\{-b\beta+a\mid a,b\in\A\}$, given by
\begin{equation}\label{eq:DL}
    D_L(x) = \begin{cases}
        -\lfloor\beta\rfloor\beta
        , & \text{if } x\in\bigl[l,\delta_{ba}\bigr]
        \text{ for } ba= \lfloor\beta\rfloor0
    ,\\
        -b\beta + a
        , & \text{if } x\in\bigl(\delta_{P(ba)},\delta_{ba}\bigr]
        \text{ for } ba\neq \lfloor\beta\rfloor0
    .\end{cases}
\end{equation}

\begin{thm}\label{t:lazybetanadruhou}
Let $\beta>1$, $\beta\notin\N$, $\A=\bigl\{0,1,\dots,\lfloor\beta\rfloor\bigr\}$.
For an $x\in I$, denote by $d_L(x)=d_{\beta^2,I,D_L}(x)$ the $\beta^2$-representation of
 $x$ obtained in accordance with Definition~\ref{JDT} by the transformation $T_L:I\to I$,
 defined by~\eqref{eq:TL}.
Then
\begin{itemize}
\item
 $d_{L}(x)$ is the lazy $\beta^2$-representation of $x$ in the alphabet $\B$;
\item
 $\psi\bigl(d_{L}(x)\bigr)$ is the lazy $(-\beta)$-representation of $x$ in the alphabet $\A$.
\end{itemize}
\end{thm}
\medskip

In \cite{ErJoKo}, the relation between the lazy and the greedy representations of numbers in a
 numeration system with a positive base $\beta \in (1,2)$ and the alphabet $\{0,1\}$
 was established.
Namely, if $x =x_1x_2x_3 \ldots$ is the greedy
$\beta$-representation of $x \in J_{\beta, \A} =
\bigl[0,\frac{1}{\beta-1}\bigr]$,
 then $(1-x_1)(1-x_2)(1-x_3)\ldots$ is the lazy $\beta$-representation of $y=\frac{ 1}{\beta-1}-x$.
Generalization of this relation for numeration systems with a base $\beta >1$ and the
 alphabet $\A=\bigl\{0,1, \dots, \lfloor\beta \rfloor\bigr\}$ was presented in \cite{DaKr}.
Comparing the formulas for $D_G$ and $D_L$, one can deduce a similar symmetry between greedy
 and lazy $(-\beta)$-representations.

\begin{prop}\label{p:symetrie}
Let $z_1z_2z_3z_4\ldots$ be the greedy $(-\beta)$-representation
of a number $z\in I$ and
 let $y_1y_2y_3y_4\ldots$ be the lazy $(-\beta)$-representation of a number $y\in I$.
Then
\[
    y_i+z_i=\lfloor\beta\rfloor
    \quad \text{for every $i\geqslant 1$}
    \qquad\iff\qquad
    y+z = -\frac{\lfloor\beta\rfloor}{\beta+1}.
\]
\end{prop}
\medskip

From now on, we shall concentrate on the properties of the greedy transformation $T_G$ only.

\begin{figure}
\centering
    \raisebox{-\height}{\includegraphics[page=9,scale=1.15]{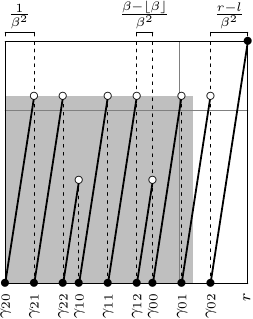}}%
\hfil
    \raisebox{-\height}{\includegraphics[page=8,scale=1.15]{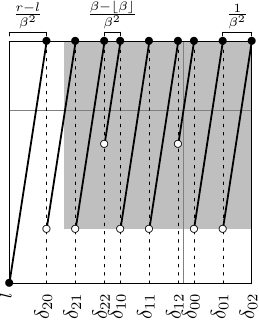}}%
\caption{The greedy and lazy transformations $T_G$ and $T_L$.
This figure corresponds to a base $-\beta\in(-3,-2)$.}
\label{graf:basebetanadruhoulazygreedy}
\end{figure}

\section{Admissibility}

The transformation $T_G$ has the following property:

\medskip

For every $x\in[l,l+1)$ one has $T_G(x)\in[l,l+1)$.
Moreover, for every $x\in I\setminus[l,l+1)$, $x\neq r$, there exists an exponent $k\in\N$
 such that $T_G^k(x)\in[l,l+1)$.

\medskip

The above fact implies that, in general, some digits from the alphabet $\B$ do not appear
 infinitely many times in the greedy $\beta^2$-representation of any number $x\in I$.

\begin{lemma}\label{l:kterecifry}
Let $d_G$ be as in Theorem~\ref{t:greedybetanadruhou}.
Denote by $\A_G$ the minimal alphabet, such that the $\beta^2$-representation $d_G(x)$ of
 every $x\in[l,l+1)$ has digits in $\A_G$.
Then
\[
    \A_G = \bigl\{-b\beta+a \mid a,b\in\A,\ b\geqslant1\bigr\} \cup \bigl\{ a
    \,\big\vert\,a\in\A,\, a < \beta\{\beta\} \bigr\},
\]
 where $\{\beta\}=\beta-\lfloor\beta\rfloor$ is the fractional part of $\beta$.
Moreover, for every $x\in I\setminus\{r\}$, the $\beta^2$-representation $d_G(x)$ has a suffix over $\A_G$.
\end{lemma}

\begin{Proof}
The fact that a given $x\in I$ satisfies $T_G^k(x)\in[l,l+1)$ for sufficiently large $k$ means
 that the greedy representation of $x$ eventually uses only such digits $-b\beta+a$ that the
 corresponding interval $[\gamma_{ba},\gamma_{S(ba)})$ has a non-empty intersection with $[l,l+1)$.
Such digits satisfy
\[
    l+1>\gamma_{ba}=\frac{b}{-\beta} + \frac{a}{(-\beta)^2}+\frac{l}{(-\beta)^2}.
\]
This is equivalent to
\begin{equation}\label{eq:ciframax}
    \beta\bigl(\beta-\lfloor\beta\rfloor\bigr) > -b\beta+a.
\end{equation}
The left-hand side of the above inequality is positive.
The right-hand side is negative whenever $b\geqslant1$.
If $b=0$ then the condition for $a$ is written as $a<\beta\{\beta\}$.
\end{Proof}

\begin{remark}\label{rem:maxdigit}
From the above lemma it follows that the maximal digit $\maxAG$ satisfies
\[
    0\leqslant \maxAG = -0\cdot\beta+\bigl\lceil\{\beta\}\beta\bigr\rceil-1.
\]
The condition $\A_G=\B$ is equivalent to the fact that the maximal
 digit $\maxB =-0\cdot\beta + \lfloor\beta\rfloor$ of the alphabet $\B$ also
 satisfies~\eqref{eq:ciframax}.
This gives the condition
\[
    \beta^2-\lfloor\beta\rfloor\beta-\lfloor\beta\rfloor>0,
\]
 which happens if and only if $\beta$ is strictly greater than the larger root
 of $x^2-\lfloor\beta\rfloor x-\lfloor\beta\rfloor$.
Therefore the condition $\A_G=\B$ is satisfied precisely by the numbers
\[
    \beta\in\bigcup_{m\in\N}\Bigl(\frac{m+\sqrt{m^2+4m}}{2},m+1\Bigr).
\]
\end{remark}

As a consequence of Lemma~\ref{l:kterecifry}, for the description of admissible digit strings
 in the greedy representations, it is reasonable to consider the transformation $T_G$ restricted
 to the interval $[l,l+1)$.
Denote for any $X=-b\beta+a\in\A_G$
\begin{equation}\label{eq:znaceniLR}
    l_X=\gamma_{ba}
    \quad\text{and}\quad
    r_X=\begin{cases}
         \gamma_{S(ba)} &\text{if } X\neq \maxAG
    ,\\
         l+1 &\text{if } X=\maxAG,
    \end{cases}
\end{equation}
 where $\gamma_{ba}=-b\beta+a$, see~\eqref{eq:gamma}.
Then we have a disjoint union
\[
    [l,l+1)=\bigcup_{X\in\A_G}[l_X,r_X)
\]
and the digit function $D_G$ restricted to $[l,l+1)$ can be written as
\[
    D_G(x)= X \quad\text{if } x\in[l_X,r_X).
\]
The mapping $D_G:[l,l+1)\to\A_G$, and the corresponding $T_G:[l,l+1)\to[l,l+1)$,
 $d_G:[l,l+1)\to\A_G^\N$ are right continuous.

\begin{figure}
\centering
    \includegraphics[page=10,scale=1.4]{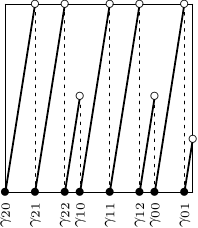}%
\qquad
    \includegraphics[page=11,scale=1.4]{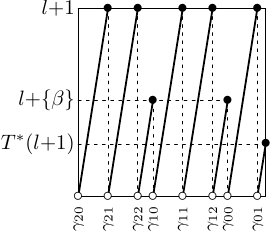}%
\caption{Right and left continuous mappings $T$ and $T^*$ defined on $[l,l+1)$ and $(l,l+1]$, respectively.
 This figure corresponds to a base $-\beta\in(-3,-2)$.}
\label{graf:TaT*}
\end{figure}

In order to formulate the result about admissible greedy representations which is derived using a
 result of~\cite{KaSt}, let us define the left continuous mapping $D^*_G:(l,l+1]\to\A_G$ by
\[
    D^*_G(x)= X \quad\text{if } x\in(l_X,r_X].
\]
The left continuous mappings $T^*_G:(l,l+1]\to(l,l+1]$, $d^*_G:(l,l+1]\to\A_G^\N$ are
 defined accordingly.
The mappings $T$ and $T^*$ are illustrated in Figure~\ref{graf:TaT*}.
Note that for every $x\in(l,l+1]$, one has
\[
    D^*_G(x)= \lim_{\varepsilon\to0+}D_G(x-\varepsilon)
    ,\quad
    T_G^*(x)=\lim_{\varepsilon\to0+} T_G(x-\varepsilon)
    ,\quad
    d_G^*(x)=\lim_{\varepsilon\to0+} d_G(x-\varepsilon)
.\]

\begin{thm}\label{t:admisgreedy}
Let $X_1X_2X_3\ldots \in \A_G^\N$. Then there exists an
$x\in[l,l+1)$ such that $d_G(x)=X_1X_2X_3\ldots$  if and only if
 for every $k\geqslant 1$
\begin{equation}\label{eq:admisgreedy}
    X_{k+1}X_{k+2}X_{k+3}\ldots \prec \begin{cases}
        d_G^*\bigl(T_G^*(l+1)\bigr)&\text{if $X_k= \maxAG$}
    ,\\
        d_G^*\bigl(l+\{\beta\}\bigr) &\text{if $X_k=-b\beta+\lfloor\beta\rfloor$, $X_k\neq \maxAG$}
    .\end{cases}
\end{equation}
\end{thm}

\begin{Proof}
Applying Theorem~2.5 from~\cite{KaSt} to the greedy transformation $T_G$, we obtain
 a necessary and sufficient condition for a string $X_1X_2X_3\ldots$ over the alphabet $\A_G$
 to be admissible, namely
\begin{equation}\label{eq:admisleft}
    d(l_{X_k})\preceq  X_kX_{k+1}X_{k+2}\ldots \prec d^*(r_{X_k})
    \qquad\text{ for all $k\geqslant 1$}
.\end{equation}

Here and in what follows, for simplicity of notation we omit in the proof the subscripts
 in $d_G=d$, $T_G=T$, $d^*_G=d^*$, $T_G^*=T^*$.

Since $d(l_{X})=Xd\bigl(T(l_{X})\bigr)$ and $d^*(r_{X})=Xd^*\bigl(T^*(r_{X})\bigr)$,
 we can rewrite conditions~\eqref{eq:admisleft} as
\begin{equation}\label{eq:admisleftlepsi}
    d\bigl(T(l_{X_k})\bigr)\preceq  X_{k+1}X_{k+2}\ldots \prec d^*\bigl(T^*(r_{X_k})\bigr)
    \qquad\text{ for all $k\geqslant 1$}
.\end{equation}

Let us study the left-hand side inequality.
From the description of $T=T_G$, it is clear that for every letter $X\in\A_G$,
 we have $T(l_X)=l$ (see also Figure~\ref{graf:TaT*}) and thus, the string on the left-hand
 side of~\eqref{eq:admisleftlepsi} reads
\[
    d\bigl(T(l_{X_k})\bigr)=d(l)=(\minAG)^\omega.
\]
Consequently, the left-hand side condition is not limiting.

\medskip
In order to study the right-hand side condition from~\eqref{eq:admisleftlepsi},
 realize that if $X\neq\maxAG$, we have
\[
    T^*(r_X) = \begin{cases}
        l+1&\text{for $X\neq -b\beta+\lfloor\beta\rfloor$}
    ,\\
        l+\{\beta\}&\text{for $X=-b\beta+\lfloor\beta\rfloor$}
    ,\end{cases}
\]
 see Figure~\ref{graf:TaT*}.
As $r_{X}=l+1$ for $X=\maxAG$, the right-hand side condition can be written as
\begin{equation}\label{eq:posun}
    X_{k+1}X_{k+2}\ldots \prec \begin{cases}
        d^*\bigl(T^*(l+1)\bigr)&\text{for $X_k= \maxAG$}
    ,\\
        d^*(l+1)&\text{for $X_k\neq -b\beta+\lfloor\beta\rfloor$, $X_k\neq \maxAG$}
    ,\\
        d^*\bigl(l+\{\beta\}\bigr) &\text{for $X_k=-b\beta+\lfloor\beta\rfloor$, $X_k\neq \maxAG$}
    .\end{cases}
\end{equation}

We have $d^*(l+1) = \maxAG \, d^*\bigl(T^*(l+1)\bigr)$.
Therefore if $X_{k+1}<\maxAG$, then the condition
\[
    X_{k+1}X_{k+2}\ldots \prec d^*(l+1) =\maxAG\, d^*\bigl(T^*(l+1)\bigr)
\]
 is satisfied.
If $X_{k+1}=\maxAG$, then the condition is equivalent to
\[
    X_{k+2}\ldots \prec d^*\bigl(T^*(l+1)\bigr)
,\] which is ensured by verifying the first of~\eqref{eq:posun}
for the shifted index $k+1$. This completes the proof. \end{Proof}

\begin{remark}\label{pozn:complementaryadmis}
Using Proposition~\ref{p:symetrie} one can derive an analogous necessary and sufficient condition
 for admissible lazy $\beta^2$-representations $X_1X_2X_3\ldots$ of numbers
 in $x\in(r-1,r]$ over the alphabet $\A_L=\lfloor\beta\rfloor-\A_G$.
\end{remark}

\section{Negative golden ratio}

Let us illustrate the previous results and their implications on the example of the
 negative base $-\beta$ where $\beta=\phi=\tfrac{1+\sqrt{5}}{2}\approx1.618$ is the golden ratio.
Real numbers representable in base $-\phi$ over the alphabet $\A=\{0,1\}$ form the
interval $J_{-\phi, \A}=I=\bigl[-1, \tfrac{1}{\phi}\bigr]=[l,r]$.

\begin{figure}
\centering
    \includegraphics[page=15,scale=0.9]{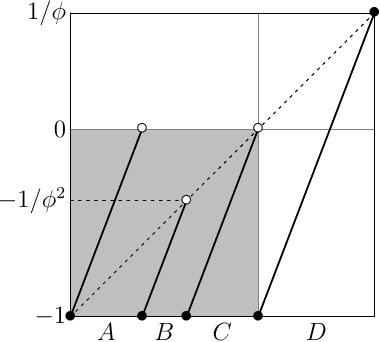}%
\hfil
    \includegraphics[page=16,scale=0.9]{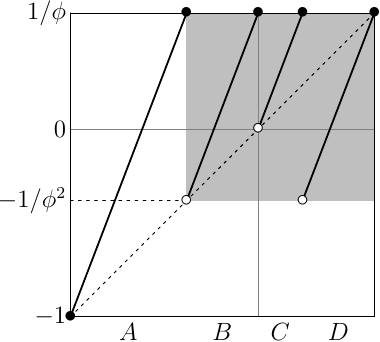}%
\caption{Greedy (left) and lazy (right) transformations $T_G$ and
$T_L$ for the base $-\phi$.}
\label{graf:LGneg}
\end{figure}

The greedy and lazy $(-\phi)$-representation can be obtained from the greedy and
 lazy $\phi^2$-representation over the alphabet $\B=\{-\phi,-\phi+1,0,1\}$,
 applying the morphism $\psi:\B^*\to\A^*$ given by
\[
    \psi(-\phi) = 10,
\quad
    \psi(-\phi+1) = 11,
\quad
    \psi(0) = 00,
\quad
    \psi(1) = 01.
\]
The greedy and lazy $\phi^2$-representations are generated by the transformation
\[
    T_G(x) = \phi^2x-D_G(x),\quad T_L(x) = \phi^2x-D_L(x),
    \quad x\in\bigl[-1, \tfrac{1}{\phi}\bigr],
\]
where the digit assigning maps $D_G$ and $D_L$ are
\[
    D_G(x) = \begin{cases}
        -\phi & \text{for } x \in {\textstyle \bigl[-1,-\frac{1}{\phi}\bigr) }
    ,\\
        -\phi+1 & \text{for } x \in {\textstyle \bigl[-\frac{1}{\phi}, -\frac{1}{\phi^2}\bigr) }
    ,\\
        0 & \text{for } x \in {\textstyle \bigl[-\frac{1}{\phi^2}, 0\bigr) }
    ,\\
        1 & \text{for } x \in {\textstyle \bigl[0,\frac{1}{\phi}\bigr] }
    ,\end{cases}
\quad
   D_L(x) = \begin{cases}
        -\phi & \text{for } x \in {\textstyle \bigl[-1,-\frac{1}{\phi^2}\bigr] }
    ,\\
        -\phi+1 & \text{for } x \in {\textstyle \bigl(-\frac{1}{\phi^2},0\bigr] }
    ,\\
        0 & \text{for } x \in {\textstyle \bigl(0,\frac{1}{\phi^3}\bigr] }
    ,\\
        1 & \text{for } x \in {\textstyle \bigl(\frac{1}{\phi^3},\frac{1}{\phi}\bigr] }
    .\end{cases}
\]
The graph of the transformations $T_G$, $T_L$ together with their restriction to the
 intervals $[l,l+1)=[-1 ,0)$, $(r-1,r]=(-\tfrac{1}{\phi^2},\tfrac{1}{\phi}]$ are drawn
 in Figure \ref{graf:LGneg}.

\begin{remark}\label{rem:G}
From the graphs of $T_G$ and $T_L$ we see that
\begin{itemize}
\item
 the greedy $\phi^2$-representation of an $x\in\bigl(-1,\frac1\phi\bigr)$
 eventually uses only digits from the alphabet $\A_G=\{-\phi,-\phi+1,0\}$, while its lazy
 $\phi^2$-representation eventually uses only digits
 $\A_L=\{-\phi+1,0,1\}$; in particular,
 \begin{itemize}
 \item
  if $x \in \bigl[0,\tfrac{1}{\phi}\bigr)$, i.\,e., outside $[l,l+1)$, then its greedy
  $(-\phi)$-representation has the form $(01)^Mz_1z_2z_3\ldots $, where $M$ is the minimal
  positive integer such that $x' = x - \sum_{k=1}^M \tfrac{1}{\phi^{2k}} <0$ and $z_1z_2z_3\ldots$
  is the greedy $(-\phi)$-representation for the number $x'$;
 \item
  if $x \in \bigl(-1,-\frac{1}{\phi^2}\bigr]$, i.\,e., outside $(r-1,r]$, then its lazy
  $(-\phi)$-representation has the form $(10)^My_1y_2y_3\ldots $, where $M$ is the minimal
  positive integer such that $x'' = x + \phi \sum_{k=1}^M \frac{1}{\phi^{2k}} > -\tfrac{1}{\phi^2}$
  and $y_1y_2y_3\ldots $ is the lazy $(-\phi)$-representation of $x''$;
 \end{itemize}
\item
 the greedy and lazy $(-\phi)$-representations of $\tfrac{1}{\phi}$ are both equal
 to $(01)^\omega$;
\item
 the greedy and lazy $(-\phi)$-representations of $-1$ are both equal
 to $(10)^\omega$;
\item
 the greedy $(-\phi)$-representation of $0$ is $01 (10)^\omega$, while its
 lazy $(-\phi)$-rep\-re\-sen\-ta\-tion is $11(01)^\omega=1(10)^\omega$.
\end{itemize}
\end{remark}

Let us now apply Theorem~\ref{t:admisgreedy} to the case $\beta=\phi$.
Denote for simplicity the digits of the alphabet $\B=\{-\phi,-\phi+1,0,1\}$ by
\[
    A=-\phi \quad < \quad B=-\phi+1 \quad < \quad C=0 \quad < \quad D=1.
\]
With this notation, we have $\A_G=\{A,B,C\}$ and $\A_L=\{B,C,D\}$.

\begin{prop}\label{p:admistau2}
A string $X_1X_2X_3\ldots$ over the alphabet $\A_G=\{A,B,C\}$ is
 the greedy $\phi^2$-representation of a number $x\in[-1,0)$ if and only if
 it does not contain a factor from the set $\{BC,B^\omega,C^\omega\}$.

A string $X_1X_2X_3\ldots$ over the alphabet $\A_L=\{B,C,D\}$ is
 the lazy $\phi^2$-rep\-re\-sen\-ta\-tion of a number
 $x\in(\tfrac{1}{\phi^2},\tfrac{1}{\phi}]$ if and only if
 it does not contain a factor from the set $\{CB,B^\omega,C^\omega\}$.
\end{prop}

\begin{Proof}
Since both $l+1=0$ and $l+\{\phi\}=-\frac1{\phi^2}$ are fixed points of the transformation $T_G^*$,
 we obtain
\[
    d_G^*\bigl(T_G^*(l+1)\bigr)=0^\omega=C^\omega
,\quad
    d_G^*\bigl(l+\{\phi\}\bigr)=(-\phi+1)^\omega=B^\omega
,\]
 and from Theorem~\ref{t:admisgreedy}, a string $X_1X_2X_3\ldots$ over the alphabet $\A_G$ is the
 greedy $\phi^2$-representation of a number $x\in[-1,0)$ if and only if for all $k\geqslant1$,
\[
    X_{k+1}X_{X+2}\ldots \prec \begin{cases}
        C^\omega & \text{ if } X_k= C
    ,\\
        B^\omega & \text{ if } X_k= B
    .\end{cases}
\]
It is not difficult to see that such condition is satisfied exactly by the strings in $\{A,B,C\}^*$
 which do not contain a factor from the set $\{BC,B^\omega,C^\omega\}$.

The statement about lazy $\phi^2$-representations is derived by
 symmetry $A \leftrightarrow D$, $B\leftrightarrow C$, (cf.\@ Proposition~\ref{p:symetrie}).
\end{Proof}

\bigskip

It follows from the facts above that only digits $B=-\phi+1$ and $C=0$ are common to both greedy
 and lazy $\phi^2$-representation restricted to the intervals $[l,l+1)$, $(r-1,r]$, respectively.
However, due to the forbidden strings $BC$ in greedy and $CB$ in lazy $\phi^2$-representation,
 no string over the alphabet $\{B,C\}$ is admissible in both cases.
Combining with Remark~\ref{rem:G}, we have the following result.

\begin{coro}
The points $x=-1$ and $x=\frac1\phi$ are the only points in $[-1,\frac1\phi]$ which have a unique
 $(-\phi)$-representation over the alphabet $\{0,1\}$.
\end{coro}

Proposition~\ref{p:admistau2} provides a combinatorial criterion for admissibility of
 representations in base $\phi^2$ in the non-integer alphabet $\A_G$.
One can also rewrite the admissibility of a digit string in base $-\phi$ using forbidden
 strings in the original alphabet $\{0,1\}$.

\begin{prop}\label{p:admistau}
A digit string $x_1x_2x_3\ldots$ over the alphabet $\{0,1\}$ is a
greedy $(-\phi)$-representation
 of some $x\!\in\![-1,0)$ if and only if
\begin{itemize}
 \item[(i)] it does not start with the prefix $1^{2k}0$, nor $0^{2k-1}1$, $k\geqslant 1$;

 \item[(ii)] it does not end with the suffix $0^\omega$ nor $1^\omega$;

 \item[(iii)] it does not contain the factor $10^{2k}1$, nor $01^{2k}0$, $k\geqslant 1$.
\end{itemize}
\end{prop}

\begin{Proof}
$(\Rightarrow)$:
Realize that grouping the digits in the greedy $(-\phi)$-representation of $x$ to pairs,
 we obtain a string over the alphabet $\{(10),(11),(00),(01)\}$.
This is justified by the fact that the greedy $(-\phi)$-representation of $x$ over
 $\{0,1\}$ is obtained from the greedy $\phi^2$-representation of $x$ over $\{A,B,C,D\}$
 applying $\psi(A)=10$, $\psi(B)=11$, $\psi(C)=00$, $\psi(D)=01$.

From Proposition~\ref{p:admistau2}, one can derive that the forbidden strings over the
 alphabet $\{(10),(11),(00),(01)\}$ are $(01)$, $(11)(00)$, $(11)^\omega$, $(00)^\omega$.
It is thus obvious that the greedy representation over the alphabet $\{0,1\}$ cannot have
 a suffix $0^\omega$, $1^\omega$, i.\,e. (ii) holds.
In order to prove (i), realize that $0^{2k-1}1 = (00)^{k-1}(01)$ contains
 the forbidden digit $(01)$.
Similarly, $1^{2k}0 = (11)^k0$ is a prefix of either  $(11)^k(01)$ or $(11)^k(00)$, where the
 former contains the forbidden digit $(01)$ and the latter the forbidden  string $(11)(00)$.
As for (iii), note that the factor $10^{2k}1$ is read in the alphabet $\{(10),(11),(00),(01)\}$
 either as $(10)(00)^{k-1}(01)$ or as a part of $(\ast1)(00)^k(1\ast)$.
In both cases it necessarily contains either $(01)$ or $(11)(00)$.
The string $01^{2k}0$ is excluded in a similar way.

$(\Leftarrow)$: We must show that a string over $\{0,1\}$
satisfying conditions (i)\,--\,(iii) does not contain
 forbidden strings of pairs $(01)$, $(11)(00)$, $(11)^\omega$, $(00)^\omega$.
Note that the brackets are always around pairs $x_{2i+1}x_{2i+2}$.
Conditions (i)\,--\,(iii) are equivalent to the fact that the
string $x_1x_2x_3\ldots$ is of the form
\begin{equation}\label{eq:label}
    0^{2k_1}1^{2k_2+1}0^{2k_3+1}1^{2k_4+1}\ldots
    ,\qquad\text{ with }k_i\geqslant0.
\end{equation}
Thus the blocks of $1$'s always start at odd-indexed positions and blocks of $0$'s
 (except the prefix) start at even-indexed positions.
This implies the impossibility of occurrence of $(01)$, $(11)(00)$.
The impossibility of $(11)^\omega$, $(00)^\omega$ is obvious.
\end{Proof}

\bigskip

Let us now compare the greedy and lazy $(-\phi)$-representations with those obtained by the
 Ito--Sadahiro algorithm, presented in Example~\ref{ItoSadahiro}.
We denote by $D_\IS$  and $T_\IS$ the corresponding transformations which now
 read $D_\IS,T_\IS:\bigl[-\frac{1}{\phi}, \frac{1}{\phi^2}\bigr)
 \mapsto\bigl[-\frac{1}{\phi}, \frac{1}{\phi^2}\bigr)$ defined
\[
    D_\IS(x) = \bigl\lfloor -\phi x  + \tfrac1\phi \bigr\rfloor
\quad\text{and}\quad
    T_\IS(x) = -\phi x - \bigl\lfloor -\phi x + \tfrac1\phi \bigr\rfloor
.\] As follows from the admissibility condition in~\cite{ItSa}, a
digit sequence $x_1x_2x_3\ldots$ over
 the alphabet $\{0,1\}$ is the Ito--Sadahiro $(-\phi)$-representation of
 some $x\in \bigl[-\frac{1}{\phi}, \frac{1}{\phi^2}\bigr)$ if and only if
\begin{equation}\label{eq:admistauIS}
    10^\omega\preceqalt x_ix_{i+1}x_{i+2}\ldots \precalt 010^\omega
    \quad\text{for all }i\geqslant1.
\end{equation}
One can easily derive that $x_1x_2x_3\ldots$ is an Ito--Sadahiro
$(-\phi)$-representation
 if and only if it does not contain any string $10^{2k+1}1$ for $k\geqslant 0$ as a substring,
 and it does not end with $010^\omega$.

\begin{figure}
\centering
    \includegraphics[page=4,scale=0.9]{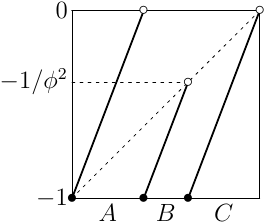}%
\hfil
    \includegraphics[page=5,scale=0.9]{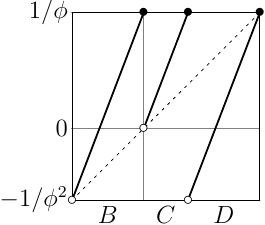}%
\hfil
    \includegraphics[page=3,scale=0.9]{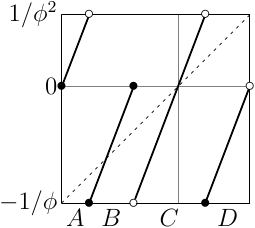}%
\caption{Transformations $T_G$ (left), $T_L$ (middle) and $T^2_\IS$ (right) in the base $\phi^2$
 that correspond to the greedy, lazy and Ito--Sadahiro representations in the base $-\phi$.}
\label{graf:maly}
\end{figure}

The example $x=-\frac12$ shows that, in general, the Ito--Sadahiro
algorithm does not assign
 to $x\in \bigl[-\frac{1}{\phi}, \frac{1}{\phi^2}\bigr)$ any extremal representation.
We have the following $(-\phi)$-representations of $-\frac12$:
\[
    \begin{array}{ll*{2}{@{\,}l}}
        \text{Ito--Sadahiro }
        & x_1x_2x_3\ldots
        & = (100)^\omega
        & = 100\,100\,100\,\ldots
    ,\\
        \text{lazy }
        & y_1y_2y_3\ldots
        & = 1(001110)^\omega
        & = 100\,111\,000\,\ldots
    ,\\
        \text{greedy }
        & z_1z_2z_3\ldots
        & = (111000)^\omega
        & = 111\,000\,111\,\ldots
    .\end{array}
\]
Clearly, $y_1y_2y_3\ldots \precalt x_1x_2x_3\ldots \precalt
z_1z_2z_3\ldots$. In fact, one can show that for no element of the
 interval $x\in \bigl[-\frac{1}{\phi}, \frac{1}{\phi^2}\bigr)$
 the Ito--Sadahiro $(-\phi)$-representation is extremal with respect to the alternate order.

\begin{proposition}\label{prop:ISnotExtremal}
The Ito--Sadahiro $(-\phi)$-representation is not extremal
 for any $x\in \bigl[-\frac{1}{\phi}, \frac{1}{\phi^2}\bigr)$.
\end{proposition}

\begin{Proof}
Let us consider an $x\in \bigl[-\frac{1}{\phi}, \frac{1}{\phi^2}\bigr)$ with
 the Ito--Sadahiro $(-\phi)$-representation $y_1y_2y_3\ldots$.
Suppose that $y_1y_2y_3\ldots$ is also a greedy
$(-\phi)$-representation of $x$. Then by Remark~\ref{rem:G}, it is
of the form
\[
    y_1y_2y_3\ldots = (01)^Mx_1x_2x_3\ldots
    ,\quad\text{for some }M\geqslant 0,
\]
where the string $x_1x_2x_3\ldots$ is the greedy representation of
some $x\in[-1,0)$. The string $x_1x_2x_3\ldots$ thus satisfies
conditions of Proposition~\ref{p:admistau},
 and, consequently, is of the form~\eqref{eq:label}.
However, this is not compatible with the admissibility condition~\eqref{eq:admistauIS}
 for Ito--Sadahiro $(-\phi)$-representations.

Suppose now that the Ito--Sadahiro $(-\phi)$-representation
$y_1y_2y_3\ldots$ is also lazy. The argumentation is similar,
noting that the forbidden strings for the
 lazy $(-\phi)$-representation can be derived from Proposition~\ref{p:admistau}
 by replacing $0\leftrightarrow1$.
\end{Proof}

\section{Unique \texorpdfstring{$(-\beta)$}{beta}-representations}
\label{sec:unique}

In~\cite{DaKa}, it is shown that for $1<\beta<2$, the set of numbers with a
 unique $(-\beta)$-representation is of Lebesgue measure zero.
The authors also show that for $\beta<\phi$, such numbers are only two.
Let us show that although the measure is always zero, the set of numbers with
 unique $(-\beta)$-representation can be uncountable.

\begin{example}\label{ex:admistribo}
Let $\beta=\mu$ be the Tribonacci constant, i.\,e., the real root
$\mu\approx1.83$ of $x^3-x^2-x-1$. Since $\lfloor\mu\rfloor=1$ and
$\mu> \phi$,  the greedy $\mu^2$-representation of numbers
 in the interval $[l,l+1)=[-\frac{\mu}{\mu^2-1},\frac1{\mu(\mu^2-1)})$ may use all possible
 digits $-\mu$, $-\mu+1$, $0$, $1$, see Lemma~\ref{l:kterecifry}.
We shall again use the notation
\[
A=-\mu, \quad B=-\mu+1, \quad C=0, \quad D=1.
\]
Our aim is to derive the admissibility condition for digit strings
$X_1X_2X_3\ldots$ over the
 alphabet $\A_G=\{A,B,C,D\}$ using Theorem~\ref{t:admisgreedy}.

\begin{figure}
\centering
    \includegraphics[page=13,scale=1.2]{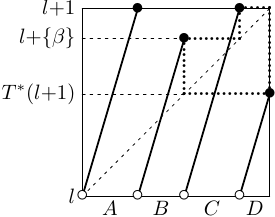}\hspace*{4em}%
\caption{The graph of the transformation $T_G^*$ for the Tribonacci constant.}
\label{graf:Tribonacci}
\end{figure}

With the notation~\eqref{eq:znaceniLR}, omitting subscripts of $T_G=T$, $d_G=d$, one
 can verify that
\begin{align*}
    l+\{\mu\} &= r_C
,\\
    T^*\bigl(l+\{\mu\}\bigr) &= T^*(r_C) = l+1 = r_D
,\\
    (T^*)^2\bigl(l+\{\mu\}\bigr) &= T^*(r_D) = r_B
,\\
    (T^*)^3\bigl(l+\{\mu\}\bigr) &= T^*(r_B) = l+\{\mu\}
,\end{align*}
 see also Figure~\ref{graf:Tribonacci}.
The strings used in the admissibility condition from Theorem~\ref{t:admisgreedy} are thus purely
 periodic and very strongly related, namely
\[
    d^*\bigl(l+\{\mu\}\bigr) = Cd^*(l+1)=CDd^*\bigl(T^*(l+1)\bigr)
    =CDBd^*\bigl(l+\{\mu\}\bigr) = (CDB)^\omega.
\]
Therefore, a digit string $X_1X_2X_3\ldots$ over the alphabet
$\A_G$ is an admissible
 greedy $\mu^2$-representation if and only if for all $k\geqslant1$,
\[
    X_{k+1}X_{k+2}X_{k+3}\ldots \prec \begin{cases}
        (CDB)^\omega & \text{ if } X_k= B
    ,\\
        (BCD)^\omega & \text{ if } X_k= D
    ,\end{cases}
\]
 which is equivalent to the requirement that the string $X_1X_2X_3\ldots$ does not contain
 a factor from $\bigl\{BD,DC,DD,(BCD)^\omega\bigr\}$.

By symmetry $A\leftrightarrow D$, $B\leftrightarrow C$, we obtain that a digit
 string $X_1X_2X_3\ldots$ over the alphabet $\A_L=\A_G$ is an admissible
 lazy $\mu^2$-representation if and only if it does not contain a factor from
 the set $\bigl\{CA,AB,AA,(CBA)^\omega\bigr\}$.

From this, we can see that any digit string over the alphabet $\{B,C\}$ is admissible both
 as the greedy and lazy $\mu^2$-representation, thus numbers represented by such a string
 have only one representation in base $-\mu$.
There are uncountably many such numbers.
\end{example}

\begin{thm}
Let $\beta>1+\sqrt{3}$, $\beta\notin\N$.
Then there exist uncountably many numbers in $J_{-\beta,\A}$ having a
 unique $(-\beta)$-representation over the alphabet $\A=\bigl\{0,1,\dots,\lfloor\beta\rfloor\bigr\}$.
\end{thm}

\begin{Proof}
We will show that there are uncountably many numbers in the interval $J_{-\beta,\A}=[l,r)$ for
 which the lazy and the greedy $(-\beta)$-representations coincide.
Necessarily, such a number $x$ (with the exception of $x=l$ and $x=r$) must lie in the
 intersection $[l,l+1)\cap(r-1,r]=(r-1,l+1)$.
The alphabet for the $\beta^2$-representation of such an $x\in(r-1,l+1)$ is $\A_G\cap\A_L$.

Recall that the requirement \eqref{eq:admisgreedy} from Theorem~\ref{t:admisgreedy} for the
 admissible greedy $\beta^2$-representations controls the string only after digits $X_k$ of
 the form $X_k=-b\beta+\lfloor\beta\rfloor$ and $X_k= \maxAG$.
By symmetry (Proposition~\ref{p:symetrie}), the requirement on lazy $\beta^2$-representations
 controls the string only after digits $X_k$ of the form $X_k=-b\beta+0$ and $X_k= \minAL$.
Realize that
\[
    \maxAG = 0\cdot\beta + a
\quad\text{ and }\quad
    \minAL = \lfloor\beta\rfloor\cdot\beta+\lfloor\beta\rfloor-a
\]
 for some $a\in\bigl\{0,1,\dots,\lfloor\beta\rfloor\bigr\}$, see Remark~\ref{rem:maxdigit}.
Therefore any digit string over the original alphabet $\A=\bigl\{0,1,\dots,\lfloor\beta\rfloor\bigr\}$
 avoiding $0$ and $\lfloor\beta\rfloor$ is allowed both as a greedy and lazy
 $(-\beta)$-representation.
If $\beta>3$, then $\A\setminus\bigl\{0,\lfloor\beta\rfloor\bigr\}$ has at least two elements
 and thus we have uncountably many numbers with unique $(-\beta)$-representation.

Let now $\beta$ satisfy $1+\sqrt3<\beta<3$.
The above argument is not sufficient, since $\A\setminus\bigl\{0,\lfloor\beta\rfloor\bigr\}=\{1\}$,
 so we need to refine our study of admissible greedy and lazy strings.
As $1+\sqrt3$ is a root of $x^2-2x-2$, by Remark~\ref{rem:maxdigit}, the alphabets $\A_G$, $\A_L$
 of greedy and lazy $\beta^2$-representations are equal to the full
 alphabet $\B=\{-b\beta+a\mid 0\leqslant a,b\leqslant2\}$.

Let us show that the string $d^*_G\bigl(l+\{\beta\}\bigr)$ starts with a
 digit $X\in\B$, $X\geqslant 0$.
For that, we have to verify that
\[
    l+\{\beta\} > l_0 = \frac{l}{(-\beta)^2},
\]
 where we use notation~\eqref{eq:znaceniLR}.
This is equivalent to $\beta^2-2\beta-2>0$, which is true precisely when $\beta>1+\sqrt3$.

From the condition~\eqref{eq:admisgreedy} we derive that any string over the alphabet
 $\{X\in\B\mid X < 0\}$ is admissible as a greedy $\beta^2$-representation.
By symmetry, any string over the alphabet $\{X\in\B\mid X> -2\beta+2\}$ is admissible as
 a lazy $\beta^2$-representation.
Altogether, every string over the alphabet
\[
    \{X\in\B\mid -2\beta+2< X < 0\} = \{-\beta,-\beta+1,-\beta+2\}
\]
is admissible both as the greedy and as the lazy $\beta^2$-representation of some
 number $x\in(r-1,l+1)$, and therefore $x$ has a unique $(-\beta)$-representation.
Obviously, there are uncountably many such numbers.
\end{Proof}

\begin{remark}
In fact, as pointed out by one of the referees, for $\beta>3$ the proof tells more
 than stated in the theorem.
The set $U$ of all numbers with $(-\beta)$-representations over the
 alphabet $\A\setminus\bigl\{0,\lfloor\beta\rfloor\bigr\}$ forms a self-similar set containing
 $\lfloor\beta\rfloor-1$ copies of itself, each time $\beta$ times smaller than the original set.
It implies that the Hausdorff dimension of the set $U$ equals to $\log_\beta\bigl(\lfloor\beta\rfloor-1\bigr)$.
The set of numbers with a unique $(-\beta)$-representation contains $U$ and therefore is of Hausdorff
 dimension at least $\log_\beta\bigl(\lfloor\beta\rfloor-1\bigr)$.
\end{remark}

\section{Conclusions}

Our main tool in this paper was to view the $(-\beta)$-representations in the
 alphabet $\A=\bigl\{0,1,\dots,\lfloor\beta\rfloor\bigr\}$ as strings of pairs of digits
 in $\A$, which amounts, in fact, to considering the alphabet $\B=-\beta\cdot \A+\A$
 and the base $\beta^2$.
Such an approach puts forward the utility of studying number systems with positive base and
 a non-integer alphabet, as was already started by Pedicini~\cite{Pedi} or Kalle and
 Steiner~\cite{KaSt}.
Obtaining new results for such systems --- for example analogous to those of de Vries and
 Komornik~\cite{Komornik} or Schmidt~\cite{Schmidt} would probably contribute also to the
 knowledge about negative base systems.

\section*{Acknowledgements}

{\small We are grateful to the referees for careful reading and
very useful suggestions.

This work was supported by the Czech Science Foundation projects
GA\v CR 201/09/0584 and 13-03538S, and by the Grant Agency of the
Czech Technical University in Prague, project
SGS11/162/OHK4/3T/14.}

\makesubmdate

\makecontacts

\end{document}